\documentclass[conference]{IEEEtran}
\IEEEoverridecommandlockouts
\usepackage{cite}
\usepackage{amsmath,amssymb,amsfonts}
\usepackage{algorithmic}
\usepackage{graphicx}
\usepackage{textcomp}
\usepackage{xcolor}

\usepackage{color}

\definecolor{red}{rgb}{1,0,0}

\newcommand{\parhead}[1]{\noindent{\bf {\em #1.}}}

\usepackage{graphicx}
\usepackage{xspace}
\usepackage{listings}

\lstdefinelanguage{java}{
morekeywords={
float, public, interface, class, static, void, extends, implements, final,
boolean, return, new, abstract, super, for, int, package, private, import,
protected, this, throw, try, finally, if, else, while, instance of, &&, ||},
	sensitive=true,
	morecomment=[l]{//},
	morecomment=[s]{/*}{*/},
	morestring=[b]",
}

\lstdefinelanguage{tracematches}
{morekeywords={
float, public, interface, class, static, void, extends, implements, final,
boolean, return, new, abstract, super, for, int, package, private, import,
protected, this, throw, try, finally, if, else, while, instance of, &&, ||,
before, around, after, returning, tracematch, sym, call, exec},
	sensitive=true,
	morecomment=[l]{//},
	morecomment=[s]{/*}{*/},
	morestring=[b]",
}

\lstdefinelanguage{JavaScript}
	{morekeywords={
	               true, false, return, new, function, var,
	               this, throw, try, finally, if, else, while, for, instanceof,
	               async, await},
	sensitive=true,
	morecomment=[l]{//},
	morecomment=[s]{/*}{*/},
	morestring=[b]",
}

\begin{document}
\title{A System Dynamics Approach to Evaluating Sludge Management Strategies in Vinasse Treatment: Cost-Benefit Analysis and Scenario Assessment}

\author{\IEEEauthorblockN{Agustin Olivares}
\IEEEauthorblockA{\textit{Escuela de Ingeniería} \\
\textit{Universidad Católica del Norte}\\
Coquimbo, Chile \\
0009-0002-6428-6805}
\and
\IEEEauthorblockN{Paul Leger}
\IEEEauthorblockA{\textit{Escuela de Ingeniería} \\
\textit{Universidad Católica del Norte}\\
Coquimbo, Chile \\
0000-0003-0969-5139}
\and
\IEEEauthorblockN{Rodrigo Poblete}
\IEEEauthorblockA{\textit{Escuela de Prevenci\'on de Riesgos y Medioambiente} \\
\textit{Universidad Católica del Norte}\\
Coquimbo, Chile \\
0000-0002-8313-1203}
}

\maketitle              
\begin{abstract}
In the Chilean local alcohol industry (pisco industry), for one liter of alcohol produced 10-15 liters of vinasse as the main wastewater of the process. To comply with industrial waste regulations, vinasse must be stored, which enables evaporation, leaving behind a residual sludge. However, treating vinasse remains an environmental and industrial challenge, having a high nutrient concentration and acidity that can degrade soil quality and harm surrounding vegetation. While previous studies have modeled sludge generation and transport in urban water systems, research on industrial wastewater, such as the alcohol industry, remains limited, affecting the search for opportunities to improve the treatment process.. This paper proposes a \textit{System Dynamics Model} (SDM) to assess the costs associated with three management strategies: natural drying of vinasse, relocation to an alternative site, and implementation of a coagulation-flocculation treatment to accelerate sludge production. This paper makes two contributions. First, we describe a pioneer SDM applicable to sludge management, which includes variables such as sludge transport, coagulant quantity, and ambient temperature to make hypothetical scenarios that affect the treatment processes of vinasse. Second, we present the expected results of the associated costs of the scenarios proposed in the model, helping decision-makers to manage vinasse. The model is calibrated with historical data provided by a company in the North of Chile, helping to improve the decision-making for vinasse treatment.
\end{abstract}



\section{Introduction}

Wastewater is a problem affecting the global population, having a water supply that is often wasted due to a lack of proper use or treatment. The main problems caused by wastewater are eutrophication and toxicity in the environment, which can lead to the death of species that inhabit those ecosystems \cite{kiely-1997}. 

Vinasse is a wastewater of the alcohol industry, generated during the distillation of ethanol from wine or other biomass feedstocks, with a high organic load, acidic pH, and significant concentrations of nutrients such as potassium, nitrogen, and phosphorus, as well as organic matter \cite{fuess-2014}. The large production volumes, with approximately 10–15 liters of vinasse per liter of ethanol as the main wastewater source from the production of alcohol \cite{christofoletti-2013}.

Wastewater is now widely perceived as a valuable resource, having treatment to recover water by removing contaminants and containing reusable components such as minerals or nutrients. \cite{villarin-2020}. The potential of wastewater to produce products from waste needs to be evaluated through wastewater treatment processes and facilities. 

For the management of wastewater, optimization, machine learning, and simulation methods are used to model their process. Optimization models have been applied to improve specific steps in the treatment train and enhance the overall operation of wastewater treatment plants \cite{wang-2021,han-2025}. Machine learning is used to integrate production and sanitation of wastewater to assess the availability to develop and optimize the location of Wastewater Treatment Plants \cite{thangarasu-2024,sakti-2023}. Similarly, Simulation models have been used to analyze wastewater treatment processes and sludge generation dynamics, and for wastewater discharges in rivers to test mitigation strategies for pollution \cite{urban-2023,pradana-2021}. Despite the advancements in modeling wastewater processes, a significant gap persists in addressing the behaviour of industrial wastewater treatment. This gap difficult to identify improvements in treatment methods for industrial applications, with challenges to find effective treatment approaches, estimating accurate costs, and developing efficient disposal methods for industrial-scale operations.

\section{Materials and methods}

A simulation model is a methodology used to represent real-world systems in a virtual environment, enabling the evaluation of hypothetical scenarios and their effects on processes. \textit{System Dynamics Modeling} (SDM) is a simulation approach that studies the relationships between variables and their mutual influences, commonly used to model the management of water supply systems \cite{urban-2023,pradana-2021}.

SDM is particularly useful for identifying relationships of variables within systems, which starts with the creation of a quantitative model using a \textit{Stock-Flow Diagram} (SFD). SFD enables the evaluation of relationships that affect the variables. SFD has been widely applied in wastewater management research, assessing the potential of artificial intelligence in industrial settings, and analyzing water resource utilization patterns for different rice genotypes \cite{yu-2024,motahari-2024}. These studies demonstrate the versatility of SFD in addressing complex water-related challenges. 

\section{Methodology}

\begin{figure}[t]
    \centering
    \includegraphics[width=0.7\linewidth]{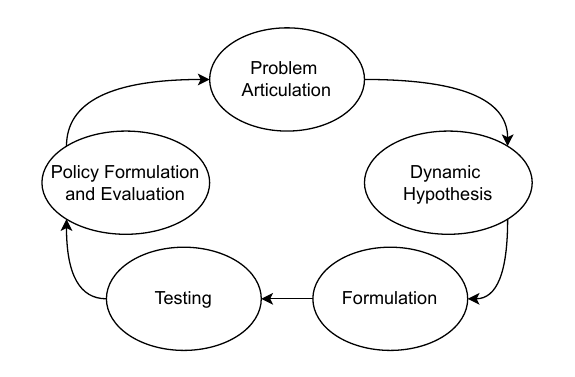}
    \caption{Process modeling of SDM}
    \label{fig:diagram-process-sd}
\end{figure}

We follow the guidelines for the modeling process of SDM mentioned in \cite{sterman-2000}. Figure~\ref{fig:diagram-process-sd} shows where the process starts with the \textit{Problem Articulation}, followed by \textit{Dynamic Hypothesis Formulation}, \textit{Simulation Model Formulation}, \textit{Testing}, \textit{Policy design and evaluation}, and \textit{Simulation Result Analysis}, as an iterative process to make improvements in the model.

The SFD presents stock variables, including vinasse accumulation, sludge accumulation, and marginal cost. SFD incorporates flows, such as vinasse rates, evaporation rates, sludge production rates, and the cost threshold for the model.

\begin{figure*}[t]
    \centering
    \includegraphics[width=0.83\linewidth]{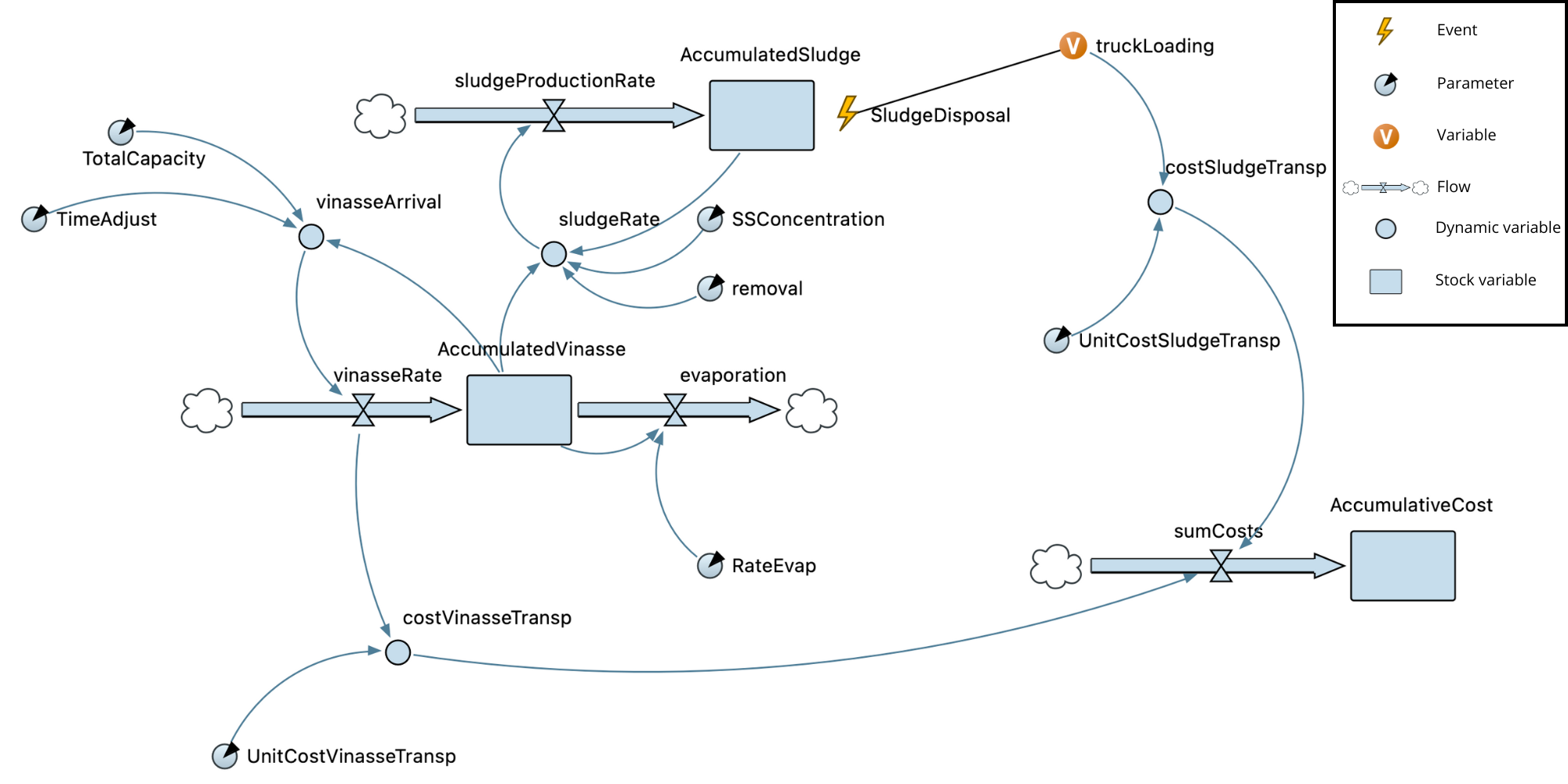}
    \caption{Stock-Flow diagram of the Vinasse plant treatment}
    \label{fig:sfd}
\end{figure*}

In Figure~\ref{fig:sfd}, the management of wastewater of the alcohol industry is modeled using an SFD. Parameters are denoted with an initial capital letter (\textit{e.g.} TotalCapacity), dynamic variables begin with a lowercase letter  (\textit{e.g.} sludgeRate), stock variables are represented as boxes (\textit{e.g.} AccumulatedSlduge), and thick arrows indicate flows  (\textit{e.g.} vinasseRate). 

The SFD is developed for quantitative analysis of the SDM by AnyLogic (8 PLE 8.9.5). The formulas used in the model rely on specific parameters, including production averages and probabilities that influence vinasse evaporation and sludge production.

In the baseline scenario, reflect the current state of the process of the vinasse treatment plant. The model begins with an empty system. The vinasse will enter the model based on the total production capacity of the plant. Upon arrival, the vinasse will be stored in a designated holding area within the plant. Simultaneously, the sludge will be generated and accumulate in the same storage pool. The vinasse contains suspended solids, which contribute to the overall composition of the stored material. This setup allows for the continuous accumulation of both vinasse and sludge. 

Future hypothetical scenarios will systematically evaluate how key operational variables influence system performance. Table~\ref{tab:scenarios} shows the scenarios to be analyzed, such as the addition of coagulants to optimize sludge generation, adjusting sludge removal efficiency by pickup frequency, and the expansion of the vinasse processing capacity plant.

The model will be tested in the vinasse treatment plan located in the Coquimbo Region of the North of Chile, for the management of byproducts from alcohol production. 

\begin{table}[t]
\centering
\caption{Model testing scenarios}
\label{tab:scenarios}
\begin{tabular}{p{0.1cm}p{2cm}p{2cm}p{2cm}}
\hline
\textbf{Nº} & \textbf{Description} & \textbf{Parameters} & \textbf{Expected Result} \\ \hline
1 & Coagulant addition & Dosage amount & Increased sludge production \\ 
& & Reduced vinasse usage \\ 
2 & Sludge transport policy & Truck capacity/frequency & More efficient removal \\ 
& & Cost analysis \\ 
3 & Plant capacity expansion & Processing capacity & Higher vinasse handling \\ 
& & Operational cost analysis \\ \hline
\end{tabular}
\end{table}

\section{Expected results}

\begin{figure}[t]
    \centering
    \includegraphics[width=0.9\linewidth]{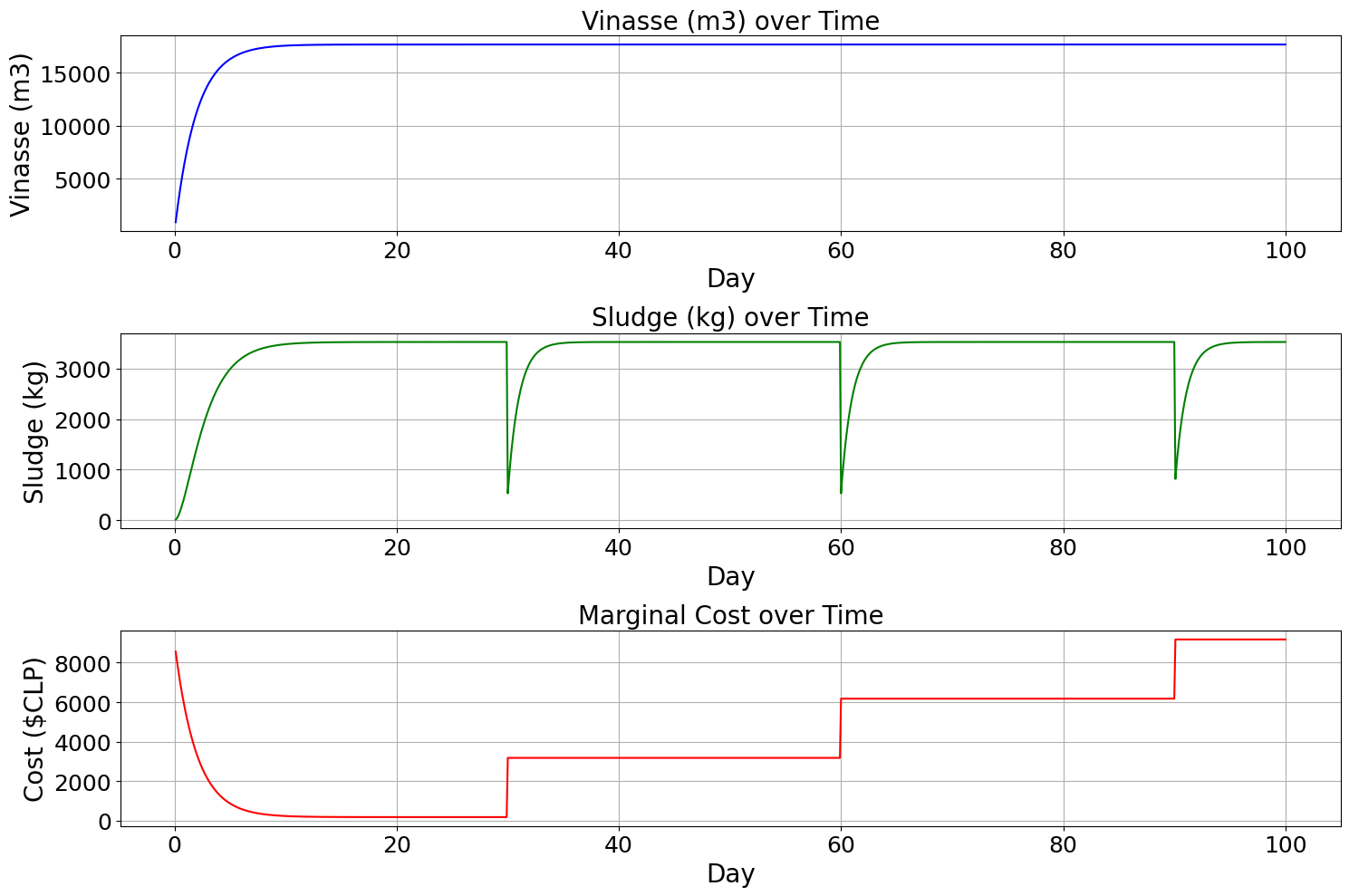}
    \caption{Results of Vinasse accumulation, Sludge accumulation, and marginal cost}
    \label{fig:results}
\end{figure}

Figure~\ref{fig:results} illustrates the initial behavior of the model, showing the values of the amount of accumulated vinasse, the amount of accumulated sludge, and the costs associated with the operation of the vinasse treatment plant. To prevent exceeding the capacity limit of 18,000 $m^3$ of vinasse of the plant, adjustments are made to the vinasse arrival rate. The Figure shows the sludge generation seasonality, with 3,000 $kg$ of sludge being extracted from the plant every 30 days. The marginal cost shows an increasing behaviour over time, maintaining a constant value until the start of sludge disposal, and subsequently increasing the transport cost of the system.

The expected results of the three proposed scenarios (coagulant addition, sludge transport policy, and plant capacity expansion) will be analyzed in comparison with the baseline model results, as detailed in the following points. 

\subsection{Coagulant addition}

The first scenario evaluates the effect of coagulant addition on sludge generation optimization. The model predicts a reduction in vinasse volume due to accelerated sludge production. However, the operational costs of coagulants may offset the savings from reduced vinasse storage requirements. Additionally, sludge extraction frequency may require adjustment to account for higher sludge yields. The results will determine the assess the feasibility of implementing a new vinasse treatment method.

\subsection{Sludge transport policy}

The next scenario examines the effects of increasing truck capacity or pickup frequency for sludge removal. This adjustment is expected to reduce sludge accumulation in the plant, thereby freeing up storage space for incoming vinasse. The model will demonstrate a more balanced system with lower peak accumulation levels. However, higher transportation costs will increase operational expenses, necessitating a cost-benefit analysis to determine the optimal removal rate. The results will support the evaluation of alternative sludge management policies to optimize storage capacity.

\subsection{Plant capacity expansion}

Expanding the vinasse processing capacity of the plant will be evaluated to assess its effect on system performance. The model will demonstrate a delayed saturation point, allowing for greater vinasse storage without immediate overflow risks. While this scenario improves operational flexibility. The results will highlight the trade-off between higher upfront costs and long-term operational advantages. The results will determine whether the capacity expansion is a cost-effective solution for improving plant performance.

\parhead{Summary} The model will analyze the combined effects of using coagulants alongside increased sludge removal efficiency, changing the transport policy, and increasing the plant capacity. The results will provide insights into the most cost-effective strategies for managing vinasse and sludge accumulation for practical solutions to the specific conditions in the vinasse treatment plant.






%
%
\bibliographystyle{splncs04}
\bibliography{references}

\end{document}